\documentclass[12pt]{article}
\usepackage{epsfig,amssymb,epic,eepic,color, graphics,amsmath}
\usepackage[english]{babel}
\input{amssym.def}

\def\nd{\noindent}

\newenvironment{demo}{\nd {\bf Proof: }}{${}$\hfill $\diamond$ \medskip}

\newtheorem{theo}{Theorem}
\newtheorem{prop}{Proposition}
\newtheorem{defi}{Definition}

\newtheorem{coro}{Corollary}

\newtheorem{rema}{Remark}

\begin{document}
\sloppy
\date{\today}
\title{Automorphism groups of rigid geometries\\ on leaf spaces of foliations}
\author{N.\,I.~Zhukova\footnote{National Research University Higher School of Economics,  
nzhukova@hse.ru}}

\maketitle

\begin{abstract} We introduce a category of rigid geometries on singular 
spaces which are leaf spaces of foliations and are considered as leaf manifolds.
We single out a special category $\mathfrak F_0$ of leaf manifolds 
containing the orbifold category as a full subcategory. Objects of $\mathfrak F_0$
may have non-Hausdorff topology unlike the orbifolds. The topology of some objects of
$\mathfrak F_0$ does not satisfy the separation axiom $T_0$. It is shown that 
for every ${\mathcal N}\in Ob(\mathfrak F_0)$ a rigid geometry $\zeta$ on $\mathcal N$
admits a desingularization. Moreover, for every such $\mathcal N$ we prove 
the existence and the uniqueness of a finite dimensional Lie group structure 
on the automorphism group $Aut(\zeta)$ of the rigid 
geometry $\zeta$ on $\mathcal{N}$.
\end{abstract}



Bibliography: \ref{last} names.

\section{Introduction and the Main results}
Singular spaces and differential geometry on them are used in
many branches of mathema\-tics and physics (see for example \cite{Bam}, \cite{Sni}, \cite{GL}). 
Orbifold, forming the full subcategory of studied in this paper category of 
leaf manifolds of foliations, used in string theory and in theory of deformation quantization.
Famous results of Thurston on the classification of closed $3$-manifolds use the classification 
of $2$-dimensional orbifolds. Orbifolds were being used by physicists in the study 
of conformal field theory, an overview of this aspect of orbifold history 
can be found in \cite{ALR}.

Different approaches to investigation of additional structures on
singular spaces of leaves of foliations are known \cite{Moer}. Grotendieck presented an approach founded
on consideration of the leaf space $M/F$ of a foliation $(M, F)$ as a topos $Sh(M/F)$ formed 
by all sheaves of $M$ which are invariant under holonomy diffeomorphisms of $(M,F)$.
Haefliger \cite{Haef} constructed and used a classifying space $B\Gamma^n$ for 
foliations of codimension $n$. Connes introduced a concept of $C^*$-algebra of 
complex valued functions with compact supports on the holonomy groupoid of a foliation 
$(M, F)$ \cite{Conn}. This $C^*$-algebra may be considered as a desingularization of 
the leaf space $M/F.$ Losik developed of some ideas of the "formal 
differential geometry" of Gel'fand \cite{Los} and applied them to the introduction new 
characteristic classes on singular leaf manifolds of foliations \cite{Los1}, \cite{Los2}. 
At present the usage of holonomy groupoids and, in particular, of \'etale 
groupoids as models of leaf spaces of foliations takes central place \cite{CM}.

As it was observed by Losik \cite{Los2}, a singular leaf space with a poor
topology may have a rich differential geometry. Our work confirms this assertion.

We investigate rigid geometrical structures on singular spaces which are 
leaf spaces of some class of smooth foliations of an arbitrary codimension 
$n$ on $(m+n)$-dimensional manifolds, where $n>0$, $m>0$. 

The rigid geometrical structures in sense of \cite{ZhR} include large classes 
of geometries such as Cartan, parabolic, conformal, projective, pseudo-Riemannian, 
Lorentzian, Riemannian, Weyl and affine connection geometries, rigid geometries 
in the sense of \cite{AG} and also $G$-structures of finite type.

In this work we introduce a concept of rigid geometries on singular 
spaces which are leaf spaces of foliations and investigate their automorphism groups.

Following Losik \cite{Los1}, we define a smooth structure on the leaf space $M/F$ 
of a foliation $(M, F)$ by an atlas (Section~\ref{ss4.1}). This smooth structure 
is called induced by $(M, F)$. Smooth leaf spaces are called by us {\it leaf manifolds}. 
The codimension of a foliation is called the dimension of the induced leaf manifold.
Leaf manifolds form a category $\mathfrak F.$ 

Further we assume that the foliations under consideration admit Ehresmann connections 
in the sense of Blumenthal and Hebda~\cite{BH}, unless otherwise specified. 
An Ehresmann connection for a foliation $(M, F)$ of codimension $n$ is an $n$-dimensional 
distribution $\mathfrak M$ transverse to $(M, F)$ which has the property of 
vertical-horizontal homotopy (we recall the exact definition in Section \ref{ss3.1}). 
An Ehresmann connection has the global differentially topological character. 

\begin{defi}\label{d1} For a given leaf manifold $\mathcal N$, a smooth foliation $(M, F)$ 
admitting an Ehresmann connection is called associated with $\mathcal N$ if 
the leaf space $M/F$ with the induced smooth structure becomes
an object of the category $\mathfrak F$ which is isomorphic to $\mathcal N$
in $\mathfrak F$.
\end{defi}

A rigid geometry on a manifold $T$ (disconnected in general) is a pair $\xi = (P(T,H),\beta)$ 
consisting of an $H$-bundle $P(T,H)$ over $T$, where $P$ is equipped with a non-degenerate 
$\mathbb R^k$-valued $1$-form $\beta$ agreed with the action of the group $H$ on $P$. 
We say that $\mathcal N\in Ob(\mathfrak F)$ has a rigid geometry $\zeta$ modelled on $\xi$ if 
there exists an associated foliation $(M, F)$ admitting $\xi$ as a transverse structure. 
Note that for a given leaf manifold $\mathcal N$, there are a lot of associated foliations
of different dimensions. We show that the definition of $\zeta$ is correct, i.e. it does 
not depend on the choice of foliation $(M, F)$ modelled on $\xi$ which is associated 
with $\mathcal N$ and we prove the following theorem.

\begin{theo}\label{t1} Let $\mathcal N$ be a leaf manifold and $(M, F)$ be an associated foliation.
Assume that $(M, F)$ is a foliation which has a transverse rigid geometry $\xi = (P(T,H),\omega)$
and admits an Ehresmann connection. Then the rigid geometry 
$\zeta = ({\mathcal R}_{\mathcal F}({\mathcal N},H),\alpha)$ on $\mathcal N$ and a structural 
Lie algebra $\mathfrak g_0 = \mathfrak g_0(\zeta)$ are defined, where ${\mathcal R}_{\mathcal F}$ 
is the leaf manifold of the lift foliation $({\mathcal R},{\mathcal F})$ for $(M, F)$ with the 
induced locally free action of the Lie group $H$ on $\mathcal R_{\mathcal F}$ such that 
${\mathcal R}_{\mathcal F}/H\cong\mathcal N$, $\alpha$ is the induced non-degenerate 
$\mathbb{R}^k$-valued $1$-form on ${\mathcal R}_{\mathcal F}$, and the Lie algebra $\mathfrak g_0$
coincides with the structural Lie algebra of $(M, F)$.
\end{theo}

The category of rigid geometries on leaf manifolds from $\mathfrak F$ is denoted by 
$\mathfrak R\mathfrak F$. The group of all automorphisms of $\xi\in\mathfrak F$
is denoted by $Aut(\xi)$ and called by the {\it automorphism group of} $\xi$.

Let $\mathfrak R\mathfrak F_0$ be the full subcategory of $\mathfrak R\mathfrak F$
objects of which have zero structural Lie algebra. Let ${\mathcal K}: \mathfrak R\mathfrak F\to
\mathfrak F$ be the covariant functor which forgets a rigid geometry. 
Put $\mathfrak F_0 = {\mathcal K}(\mathfrak R\mathfrak F_0)$ and note that 
$\mathfrak F_0$ is a full subcategory of $\mathfrak F$.

Emphasize that any $n$-dimensional orbifold belongs to $Ob(\mathfrak F_0),$ and 
$\mathfrak F_0$ is a great expansion of the orbifold category. In particular, 
every leaf manifold ${\mathcal N}\in Ob(\mathfrak F)$ admitting a rigid geometry
and satisfying the separation axiom $T_0$, belongs to  $Ob(\mathfrak F_0).$
Moreover, there are ${\mathcal N}\in Ob(\mathfrak F_0)$ which do not 
satisfy the separation axiom $T_0$.

`````````````````

We prove the following two theorems.

\begin{theo}\label{t2} Let $\zeta\in Ob(\mathfrak R\mathfrak F_{0})$ be a rigid geometry on 
$n$-dimensional leaf manifold $\mathcal N\in Ob(\mathfrak F_{0})$. 
Let $(M, F)$ be an associated foliation, $({\mathcal R},{\mathcal F})$ be its 
lifted foliation with the projection $\pi:{\cal R}\to M$ of the $H$-bundle and 
$\omega$ be the induced $\mathbb R^k$-valued $1$-form on $\mathcal R$. Then the rigid
geometry $\zeta = ({\mathcal R}_{\mathcal F}({\mathcal N},H),\alpha)$ on $\mathcal N$ 
has the following properties:
 
\begin{itemize} \item[(i)] 1) the leaf manifold ${\mathcal R}_{\mathcal F}\cong W$ is a smooth manifold 
with a smooth locally free action of the structural Lie group $H$ such that ${\mathcal N}$ 
is the orbit space $W/H$, 2) the canonical projections 
$\pi_b:{\mathcal R}\to{\mathcal R}/{\mathcal F}\cong W $,
$\pi_{\mathcal F}:W\to W/H\cong\mathcal N$ and $r:M\to M/F\cong W/H$  satisfy the equality 
$\pi_{\mathcal F}\circ\pi_b = r\circ\pi,$ 3) $\alpha$ is $\mathbb R^k$-valued non-degenerate
$1$-form on $W$ such that $\pi_{b}^{*}\alpha = \omega$,
where $\pi_{b}^{*}$ is the codifferential of $\pi_{b}$;
\item[(ii)] the automorphism group $Aut(\zeta)$ of $\zeta$
admits a structure of a finite dimension Lie group, and
its Lie group structure is defined uniquely;
\item[(iii)] the dimension of $Aut(\zeta)$ satisfies the inequality
$$\dim Aut(\zeta)\leq\dim W = k,$$
and $k = n+s$, where $s$ is the dimension of the Lie group $H$.
\end{itemize}
\end{theo}

Thus a rigid geometry on every ${\mathcal N}\in Ob(\mathfrak F_0)$ admits
desingularization indicated in Statement $(i)$ of Theorem~\ref{t2}.

Using Theorem \ref{t2} we prove the following. 
\begin{theo}\label{t3}  Let $\mathcal N$ be a leaf manifold.
Suppose that the underlying topological space of $\mathcal N$ satisfies 
the separation axiom $T_0$ and $\mathcal N$ admits a rigid geometry $\zeta$.
Then:

1) the pair  $({\mathcal N},\zeta)$ satisfies Theorem~\ref{t2};

2) there exists an open dense subset $\mathcal N_0$ of $\mathcal N$ such that
$\mathcal N_0$ with induced smooth structure is isomorphic in the category $\mathfrak F$
to an $n$-dimensional manifold, which is not necessarily connected and not necessarily Hausdorff.
\end{theo}

It is well known that for any smooth orbifold $\mathcal N$ 
there exists a Riemannian foliation $(M, F)$ with an Ehresmann connection for which it is the 
leaf space (see, for example \cite{ZhK}). This fact implies that $\mathcal N$ is a leaf 
manifold having $(M, F)$ as the associated foliation, and $(M, F)$ is a proper foliation with 
only closed leaves. Therefore for any rigid geometry $\zeta$ on an orbifold $\mathcal N$
it is necessarily $\mathfrak g_0 = 0,$ hence orbifolds form a full subcategory of $\mathfrak F_0$.

The application Theorem~\ref{t2} gives the following two statements.

\begin{theo}\label{t4}  Let $\mathcal N$ be an $n$-dimensional orbifold equipped with
a rigid geometry $\zeta = ({\mathcal R}_{\mathcal F}({\mathcal N},H),\alpha)$. Then 
the automorphism group $Aut(\zeta)$ of $\zeta$ admits a structure of a finite 
dimension Lie group, and its Lie group structure is defined uniquely, the dimension of 
$Aut(\zeta)$ satisfies the inequality 
$$\dim Aut(\zeta)\leq\dim W = n+s,$$
where $s$ is the dimension of the structural Lie group $H$ of $\zeta$.
\end{theo}

\begin{coro} (\cite[Theorem 1]{BZ}) Let $Aut(\zeta)$ be the automorphism group of a $G$-structure 
$\zeta$ of finite type and order $m$ on a smooth $n$-dimensional orbifold $\mathcal N$. Then the 
group $Aut(\zeta)$ admits a unique topology and a unique smooth structure that makes it into a 
Lie group, and the dimension of $Aut(\zeta)$ satisfies the inequality
$$\dim Aut(\zeta)\leq\dim W = n + \dim\mathfrak g + \dim\mathfrak g_1 + ... + \dim\mathfrak g_{m-1},$$
where $g_i$ is the $i$-th prolongation of the Lie algebra $\mathfrak g$ of the group $G.$
\end{coro}
The classical theorems of Myers and Steenrod, Nomizu, Hano and Morimoto, Ehresmann on 
the existence of a Lie group structure in the full automorphism groups of Riemannian, 
affine connection geometries and of a finite type structure on manifolds, 
respectively, follow from Theorem \ref{t4}.

\medskip{\noindent\bf Assumptions\,} Throughout this paper
we assume for simplicity that all manifolds and  maps are smooth
of the class $C^r$, $r\geq 1$, and $r$ is large enough which is necessary 
for a suitable rigid geometry. All neighborhoods are
assumed to be open and all manifolds are assumed to be Hausdorff
unless otherwise specified.

\medskip{\noindent\bf Notations\,} Let
$\mathfrak X(T)$ denote the module of smooth vector fields over the ring of 
smooth functions on a manifold $T.$ If $\mathfrak M$ is a smooth distribution on $M$
and $f: K\to M$ is a submersion, then let $f^*\mathfrak M$ be the
distribution on the  manifold $K$ such that $(f^*\mathfrak M)_z=
\{X\in T_zK\,|\, f_{*z}(X)\in\mathfrak M_{f(z)}\}$, where $z\in
K$. Let $\mathfrak X_{\mathfrak M}(M)=\{X\in\mathfrak X(M)\mid
X_u\in {\mathfrak M}_u\quad\forall u\in M\}$. Let $id_M$ be
the identity mapping of a manifold $M$. Denote by $\mathcal Fol$
the foliation category in which morphisms are smooth maps transforming leaves to leaves.

The symbol  $\cong$ will denote the isomorphism of objects in the
corresponding category.

{\it Acknowledgements.} I express my gratitude to Anton Galaev 
who drew my attention to the works of M.~V.~Losik.

The publication was prepared within the framework of the Academic Fund Program 
at the National Research University Higher School of Economics (HSE) in 2016--2018 
(grant No 16-01-0010) and by the Russian Academic Excellence Project "5-100".

\section{Rigid geometries}
\subsection{Rigid structures} 
A manifold that admits an $e$-structure is called parallelizable. In other
words, a parallelizable manifold is a pair $(P,\beta),$ where $P$
is a $k$-dimensional smooth manifold and $\beta$ is a smooth non-degenerate
$\mathbb{R}^k$-valued 1-form $\beta$ on $P,$ i.e.,
$\beta_u\colon T_uP\to{\Bbb R}^k$ is an isomorphism of the vector
spaces for each $u\in P.$ 

Denote by $P(T,H)$ a principal $H$-bundle with the projection
$p\colon P\to T$. Suppose that the action of $H$ on $P$ is a right
action and $R_a$ is the diffeomorphism of $P$ corresponding to an
element $a\in H.$

Two principal bundles $P(T,H)$ and $\tilde P(\tilde T,\tilde H)$
are called {\it isomorphic} if $H=\tilde H$ and there exists a
diffeomorphism $\Gamma\colon P\to\tilde P$ such that $\Gamma\circ
R_a=\tilde R_a\circ\Gamma\,\,$ $\forall a\in H,$ where $\tilde R_a$
is the transformation of $\tilde P$ defined by an element
$a\in H.$

\begin{defi}\label{d2} Let $P(T,H)$ be a principal
$H$-bundle and $(P,\beta)$ be a parallelizable manifold
satisfying the following condition:
\begin{itemize}
\item[(S)] there is an inclusion $\mathfrak h\subset\mathbb{R}^k$ of
the vector space of the Lie algebra $\mathfrak h$ of the Lie group $H$
into the vector space $\mathbb{R}^k$ such that
$\beta(A^*)=A\,\, \forall A\in\mathfrak h,$ where $A^*$ is the
fundamental vector field on $P$ corresponding to $A.$
\end{itemize}
Such $\xi=(P(T,H),\beta)$ is called a {\it rigid structure} on
the manifold $T.$ A pair $(T,\xi)$ is called a {\it rigid
geometry}. The Lie group $H$ is called the structural Lie group of $\xi$.
\end{defi}

\begin{defi}\label{d3} Let $\xi=(P(T,H),\beta)$ and 
$\xi'=(P'(T', H),\beta')$ be two
rigid structures with the projections $p:P\to T$ and  $p':P'\to T'.$
An isomorphism $\Gamma\colon P\to P'$ of the 
$H$-bundles $P(T,H)$ and $P'(T',H)$
satisfying the equality $\Gamma^*\beta'=\beta$ is called an
isomorphism of the rigid structures $\xi$ and $\xi'.$

Such isomorphism $\Gamma$ defines a map $\gamma\colon T\to T'$ 
satisfying the equality $p'\circ\Gamma=\gamma\circ p,$ and $\gamma$ is 
a diffeomorphism of $T$ onto $T'.$ The projection $\gamma$ is called an {\it
isomorphism} of the rigid geometries $(T,\xi)$ and $(T',\xi').$
\end{defi}

\subsection{Induced rigid geometries}
Let $\xi=(P(T,H),\beta)$ be a rigid structure on a manifold $T$
with the projection $p\colon P\to T.$ Let $V$ be an arbitrary open
subset of the manifold $T,$ let $P_V:=p^{-1}(V)$ and
$\beta_V:=\beta|_{P_V}.$ Then $\xi_V:=(P_V(V,H),\beta_V)$ is
also a rigid structure.

\begin{defi}\label{d4} The pair $(V,\xi_V)$ defined above is
called the {\it induced rigid geometry on the open subset} $V$ of
$T.$
\end{defi}

\subsection {Effectiveness of rigid geometries}

Let $Aut(\xi)$ be the automorphism group of a rigid
structure $\xi=(P(T,H),\beta).$ It is a Lie group as a closed
subgroup of the Lie automorphism group $Aut(P,\beta)$ of the parallelizable 
manifold $(P,\beta).$ Denote by $Aut(T,\xi)$ the group 
of all automorphisms of the geometry $(T,\xi),$ i.e., 
$Aut(T,\xi):=\{\gamma\in\mathrm{Diff}(T)\mid\exists\,\Gamma\in{Aut}(\xi): 
p\circ\Gamma=\gamma\circ p\}.$ Consider the group epimorphism 
$\chi\colon Aut(\xi)\to Aut(T,\xi)\colon\Gamma\mapsto\gamma,$ where
$\gamma$ is the projection of $\Gamma$ with respect to $p\colon
P\to T.$

\begin{defi}\label{d5} Let $\xi=(P(T,H),\beta)$ be a rigid
structure on a manifold $T$ and let $p\colon P\to T$ be the projection.
The group $\mathrm{Gauge}(\xi):=\{\Gamma\in{Aut}(\xi)\mid
p\circ\Gamma=p\}$ is called the {\it group of gauge transformations
of the rigid structure} $\xi.$
\end{defi}

Remark that $\mathrm{Gauge}(\xi)$ is a closed normal Lie subgroup
of the Lie group ${Aut}(\xi)$ as the kernel of the group epimorphism
$\chi\colon{Aut}(\xi)\to{Aut}(T,\xi).$

\begin{defi}\label{d6} A rigid structure
$\xi=(P(T,H),\omega)$ is called {\it effective} if for any
open connected subset $V$ in $T$ the induced rigid structure
$\xi_V=(P_V(V,H),\beta_V)$ has the trivial group of gauge
transformations, i.e.,
$\mathrm{Gauge}(\xi_V)=\{\mathrm{id}_{P_V}\}.$ A rigid geometry $(T,\xi)$ is said 
{\it to be effective} if $\xi$ is an effective rigid structure.
\end{defi}

We want to emphasize that effective Cartan geometries (\cite{Shar}, \cite {C-S}), 
$G$-structures of finite type and rigid structures in the sense of \cite{AG} are examples of 
effective rigid geometries.

\subsection {Pseudogroup of local automorphisms} 
Let $(T,\xi)$ be a rigid geometry, and the topological space of $T$
may be disconnected. For arbitrary open subsets $V,$ $V'\subset T$ an
isomorphism $V\to V'$ of the induced rigid geometries $(V,\xi_V)$
and $(V',\xi_{V'})$ is called a {\it local automorphism} of
$(T,\xi).$ The family ${\cal H}$ of all local automorphisms of a
rigid geometry $(T,\xi)$ forms a pseudogroup of local
automorphisms. Denote it by ${\cal H}={\cal H}(T,\xi).$ Recall
that a pseudogroup ${\cal H}$ of local diffeomorphisms of manifold
$T$ is called {\it quasi-analytic} if the existence of an open
subset $V\subset T$ and an element $\gamma\in{\cal H}$ such that
$\gamma|_V=\mathrm{id}_V$ implies that
$\gamma|_{D(\gamma)}=\mathrm{id}_{D(\gamma)}$ in the entire
(connected) domain ${D(\gamma)}$ on which $\gamma$ is defined.

The following statement is important in the future.

\begin{prop}\label{p1} A pseudogroup ${\cal H}={\cal
H}(T,\xi)$ of all local automorphisms of an effective rigid
geometry $(T,\xi)$ is quasi-analytic.
\end{prop}
\begin{demo} Let $\gamma\in{\cal H}$ be defined on an open connected subset $D(\gamma)$
of $T$. Assume that there exists an open subset $V\subset D(\gamma)$ such that 
$\gamma|_V=\mathrm{id}_V$. Since $\gamma$ is the projection of an local automorphism
$\Gamma: P_{D(\gamma)}\to P_{D(\gamma)}$ of the rigid geometry $\xi$, it is 
necessary $\Gamma\in\mathrm{Gauge}(\xi_V)$. According to  Definition~\ref{d6} of effectiveness
of the rigid geometry $\xi$ we have $\Gamma_{P_V} = id_{P_V}$. As $P_V\subset P_{D(\gamma)}$
and $D(\gamma)$ is connected, then each connected component $P_{D(\gamma)}^c$ of $P_{D(\gamma)}$
contains some connected component $P_{V}^c$ of $P_{V}$. Note that $\Gamma$ preserves $P_{D(\gamma)}^c$.
It is well known that every automorphism of a connected parallelizable manifold is uniquely determined 
by its value at a single point. Since $\Gamma(w) = w$ for $w\in P_V^c$, then $\Gamma = id_{P_{D(\gamma)}^c}$,
where $P_{D(\gamma)}^c$ is an arbitrary connected component of $P_{D(\gamma)}$. Hence
$\Gamma = id_{P_{D(\gamma)}}$. This implies $\gamma = id_{D(\gamma)}$.
\end{demo}

\section{Foliations with transverse rigid geometries. \\Foliated bundles}\label{ss3.1}

\subsection{Ehresmann connections for foliations}\label{ss3.1}

The notion of an Ehresmann connection for foliations was introduced 
by Blumenthal and Hebda~\cite{BH}. We use the terminology from \cite{Min}. 
Let $(M, F)$ be a smooth  foliation of codimension $n\geq 1$ and  $\mathfrak{M}$ be an
 $n$-dimensional  transversal distribution on $M$. All maps and curves considered 
here are assumed to be piecewise smooth. The curves in the leaves of the 
foliation are called vertical; the distribution $\mathfrak{M}$ and its  
integral curves are called horizontal.

A map $H:I_1\times I_2\to M,$ where
 $I_1=I_2=[0,1]$, is called a vertical-horizontal homotopy if
for each fixed  $t\in I_2$, the curve $H_{|I_1\times \{t\}}$ is
horizontal, and for each fixed  $s\in I_1$, the curve
$H_{|\{s\}\times I_2}$ is vertical. The pair of curves
$(H_{|I_1\times \{0\}},H_{|\{0\}\times I_2})$ is called the base
of $H$.

A pair of curves  $(\sigma,h)$ with a common starting point
$\sigma(0) = h(0)$, where $\sigma:I_1\rightarrow M$
 is a horizontal curve, and $h:I_2\rightarrow M$ is a vertical curve, is called admissible.
If for  each admissible pair of curves $(\sigma,h)$ there exists a
vertical-horizontal homotopy with the base $(\sigma,h)$, then the
distribution  $\mathfrak M$ is called an {\it Ehresmann connection for
the  foliation} $(M, F)$. Note that there exists at most one
vertical-horizontal homotopy with a given base. 

\subsection{Foliations with transverse rigid geometries}

Let $(T,\xi)$ be a rigid geometry on an $n$-dimensional manifold $T,$
and the topological space of $T$ may be disconnected.
A foliation $(M, F)$ of codimension $n$ on an $(m+n)$-dimensional
manifold $M$ has a transverse rigid geometry $(T,\xi)$ 
if $(M, F)$ is 
defined by a cocycle $\eta=\{U_i,f_i,\{\gamma_{ij}\}\}_{i,j\in J}$ 
modelled on $(T,\xi),$ i.e.,
\begin{itemize} \item[1)] $\{U_i\,|\, i\in J\}$ is an open covering of $M;$
\item[2)] $f_i\colon U_i\to T$ are submersions with connected
fibres; \item[3)] $\gamma_{ij}\circ f_j=f_i$ on $U_i\cap U_j,$ where 
$$\gamma_{ij}: (f_{j}(U_i\cap U_j), \xi_{f_{j}(U_i\cap U_j)})\to
(f_{i}(U_i\cap U_j), \xi_{f_{i}(U_i\cap U_j)})$$ 
is a local automorphism of $(T,\xi).$ 
\end{itemize}

Without loss of generality, we will suppose that $T=\cup_{i\in J}f_i(U_i)$ 
and the family $\{(U_i,f_i)\}_{i\in J}$ is maximal
as it is generally used in the manifold theory.

For short $(M, F)$ is referred to as a foliation with TRG (i.e. with a transverse rigid geometry).

Recall that the pseudogroup generated by local diffeomorphisms $\gamma_{ij}$,
$i,j\in J$, is referred to as the {\it holonomy pseudogroup} of $(M, F)$.
It is denoted by $\mathcal H = \mathcal H(M, F).$ 

\begin{defi}\label{d7} The cocycle $\eta$ modelled on $(T,\xi)$ is
said to be an {\it $(T,\xi)$-cocycle}.
It is said also that $(M, F)$ is modelled on the rigid geometry $(T,\xi)$.
\end{defi}

Note that an $e$-foliation (or a transversally parallelizable foliation) is a foliation 
admitting a transverse rigid geometry with the trivial structure Lie group $H$, i.e.
$H = \{e\}$.

\subsection{The lifted foliation}\label{ss3.3}

We use the construction of the lifted foliation $(\cal R,
\cal F)$ for a foliation $(M, F)$ with TRG from \cite{ZhR}. It
generalizes a similar construction for an effective Cartan foliation \cite{Min}
and for a Riemannian foliation \cite{Mo}. For a given foliation $(M, F)$
with TRG one may construct a principle $H$-bundle ${\mathcal R}(M,H)$  (called
a foliated bundle) with the projection $\pi:{\mathcal R}\to M$, an
$H$-invariant transversally parallelizable foliation $({\mathcal
R},{\mathcal F})$ such that  $\pi$ is a morphism of $({\mathcal R},{\mathcal
F})$ onto $(M, F)$ in the category of foliations $\mathcal Fol$. Moreover, there
exists a $\mathbb R^k$-valued $1$-form $\omega$ on ${\mathcal R}$ having
the following properties:

(i) $\omega(A^*)=A$ for any $A\in \mathfrak h,$ where $A^*$ is the
fundamental vector field corresponding to $A$;

(ii) for any $u\in {\mathcal R}$, the map $\omega_u: T_u{\mathcal R}\to
\mathbb R^k$ is surjective with the kernel $\ker \omega = T{\mathcal
F}$, where $T{\mathcal F}$ is the tangent distribution to the
foliation $(\mathcal R, \mathcal F)$;

(iii) the Lie derivative $L_X\omega$ is zero for any vector field $X$ 
tangent to the leaves of $(\cal R, \cal F).$

The foliation $(\mathcal R, \mathcal F)$ is called the {\it lifted
foliation}. 

The restriction $\pi|_{\mathcal L}: {\mathcal L}\to L$ of $\pi$ to a leaf 
$\mathcal L$ of $(\mathcal R, \mathcal F)$ is a regular covering map onto 
the corresponding leaf $L$ of $(M, F)$, and the group of deck transformation 
of $\pi|_{\mathcal L}$ is isomorphic to the germ holonomy group of $L$
which is usually used in the foliation theory. 

If $\mathcal R$ is disconnected, then we consider a connected component of 
$\mathcal R$ and denote it also $\mathcal R$. Without loss of generality 
we assume that the same Lie group $H$ acts on $\mathcal R$ and $P$.

Let $(M, F)$ be defined by a $(T,\xi)$-cocycle 
$\eta=\{U_i,f_i,\{\gamma_{ij}\}\}_{i,j\in J}.$ Effectiveness of 
$\xi$ guarantees the existence of a unique isomorphism $\Gamma_{ij}$ 
of the induced rigid structures $\xi_{f_j(U_i\cap U_j)}$ and $\xi_{f_i(U_i\cap
U_j)}$, whose projection  coincides with $\gamma_{ij}$. Hence, in
the case $U_i\cap U_j\cap U_k\neq\emptyset,$ the equality
$\gamma_{ij}\circ\gamma_{jk}=\gamma_{ik}$ implies the equality
$\Gamma_{ij}\circ\Gamma_{jk}=\Gamma_{ik}.$
The following two equalities are direct corollaries of the previous equality
and the effectiveness of $\eta$: 
$\Gamma_{ii}=\mathrm{id}_{P_i}$ and $\Gamma_{ij}=(\Gamma_{ji})^{-1}.$

Remark that the holonomy pseudogroup of $(\mathcal R, \mathcal F)$ is generated by 
$\Gamma_{ij}$, $i,j\in J$.

\subsection{The structural Lie algebra of an $e$-foliation  
with an Ehresmann connection}\label{ss3.4}
 
At first we prove the following theorem.
\begin{theo}\label{t5} Let $(M, F)$ be an $e$-foliation with an Ehresmann
connection on a connected manifold $M$. Then the closures of its leaves are
fibers of a locally trivial fibration $\pi_b: M\to W$ over the manifold $W$.
On every fiber $\overline{L}$, $L\in F$, of this fibration the induced foliation
$(\overline{L},F|_{\overline{L}})$ is a Lie foliation with dense leaves.
The structural Lie algebra $\mathfrak g_0$ of the Lie foliation $(\overline{L},F|_{\overline{L}})$ 
does not depend on the choice of $L\in F$ and $\mathfrak g_0$  is called the structural Lie 
algebra of $(M, F)$. 
\end{theo}
\begin{demo} Since an $e$-foliation $(M, F)$ is a Riemannian one and
it admits an Ehresmann connection, according to \cite[Proposition 2]{ZhG} 
the holonomy group ${\mathcal H} = {\mathcal H}(M, F)$ is complete. Therefore we may apply 
the results of Salem \cite{Sal}. According to \cite{Sal}, the closure $\overline{L}$ of a leaf $L$
is a smooth manifold and the induced foliation $(\overline{L},F|_{\overline{L}})$ is 
a Lie foliation with dense leaves. Let $\mathfrak M$ be an Ehresmann connection
for $(M, F)$, then $\mathfrak N:= \mathfrak M\cap\overline{L}$ is an Ehresmann 
connection for $(\overline{L},F|_{\overline{L}})$. By 
\cite[Proposition 2]{ZhG}, the pseudogroup of $(\overline{L},F|_{\overline{L}})$
is complete. Therefore the structural Lie algebra $\mathfrak g_0 = \mathfrak g_0(\overline{L},F|_{\overline{L}})$  of 
the Lie foliation $(\overline{L},F|_{\overline{L}})$ is defined \cite{Sal}. Since the foliation $(M, F)$ admits 
an Ehresmann connection, the automorphism group of $(M, F)$ in the foliation category $\mathcal Fol$ 
acts transitively on the set of its leaves. This means that for every leaves $L$ and $L'$ there is an
automorpism $f: M\to M$ of $(M, F)$ such that $f(L) = L'$. The diffeomorphism $f$ has the property
$f(\overline{L}) = \overline{L'}$. Therefore, $f|_{\overline{L}}$ is an isomorphism of the Lie foliations 
$(\overline{L}, F_{\overline{L}})$ and $(\overline{L'}, F_{\overline{L'}})$ in $\mathcal Fol$.
It is well known that the structural Lie algebra of a Lie foliation with dense leaves and with 
the complete holonomy pseudogroup is an invariant in the category $\mathcal Fol$ \cite{Hae}, \cite{Sal}. 
Therefore $\mathfrak g_0(\overline{L}, F_{\overline{L}}) = \mathfrak g_0(\overline{L'}, F_{\overline{L'}})$
and the definition $\mathfrak g_0(M, F): = \mathfrak g_0(\overline{L}, F_{\overline{L}})$ is
correct.

Observe that the foliation $(M, \overline{F})$ formed by closures of leaves of $(M, F)$
is a regular Riemannian foliation with an induced Ehresmann connection, and all its leaves 
are closed. By analogy with proof of \cite[Theorem 4.2']{Mo} we show that a leaf 
$\overline{L}\in\overline{F}$ has a saturated neighborhood $\mathcal U$ such that
$({\mathcal U}, \overline{F}_{\mathcal U})$ is $e$-foliation. This implies that leaves of 
$(M, \overline{F})$ are fibers of a locally trivial fibration which is denoted by  $\pi_b: M\to W$.
\end{demo}

For a complete $e$-foliations a similar theorem is proved in \cite[Theo\-rem~4.2']{Mo}.

\subsection{The structural Lie algebra of a foliation with TRG admitting an Ehresmann connection}
We use notations introdused in Section~\ref{ss3.3}.
Let $(M, F)$ be a foliation with TRG having an Ehresmann conection $\mathfrak M$.
Let $(\mathcal R, \mathcal F)$ be the lifted $e$-foliation and $\pi:{\mathcal R}\to M$
be the projection of $H$-bundle ${\mathcal R}(M,H)$. Observe that
the distribution $\widetilde{\mathfrak M} = \pi^*\mathfrak M$ is an Ehresmann connection for  
the foliation $(\mathcal R, \mathcal F)$. Applying Theorem~\ref{t5}
to $e$-foliation $(\mathcal R, \mathcal F)$ admitting an Ehresmann connection
we obtain the following statement.

\begin{theo}\label{t6} Let $(M, F)$ be a foliation
with TRG admitting an Ehresmann connection and let $({\cal R},{\cal F})$ 
be its lifted $e$-foliation. Then:
\begin{itemize} \item[(i)] the closures of the leaves of the foliation
$\cal F$ are fibers of a certain locally trivial fibration
$\pi_b\colon{\cal R}\to W;$ \item[(ii)] the foliation
$(\overline{\cal L},{\cal F}|_{\overline{\cal L}})$ induced on the
closure $\overline{\cal L}$ is a Lie foliation with dense leaves
with the structure Lie algebra $\frak g_0$, that is the same for
any $\cal L\in\cal F.$
\end{itemize} 
\end{theo} 
According to Theorem~\ref{t6} the following definition is correct.

\begin{defi}\label{d8} The structural Lie algebra $\frak g_0$
of the Lie foliation $(\overline{\cal L},{\cal F}|_{\overline{\cal
L}})$ is called {\it the structural Lie algebra} of the foliation $(M,F)$ 
with TRG admitting an Ehresmann connection
and is denoted by $\frak g_0=\frak g_0(M,F).$
\end{defi}

\begin{rema}\label{r1} If $(M, F)$ is a Riemannian foliation on a
compact manifold, this notion coincides with the notion of a
structural Lie algebra in the sense of Molino~\cite{Mo}.
\end{rema}

\begin{defi}\label{d9} The fibration $\pi_b\colon{\cal
R}\to W$ satisfying Theorem~\ref{t5} is called a {\it basic fibration}
for $(M, F).$
\end{defi}
\begin{rema}\label{r2} Under stronger conditions of completeness of $(M, F)$ 
a similar theorem is obtained in \cite[Theorem 2]{ZhR}.
The advantage of Theorem~\ref{t6} in comparing with \cite[Theorem 2]{ZhR} 
is also that the condition of the existence of an Ehresmann connection for $(M, F)$ is defined on $M$ 
in contrast to the completeness of $(M, F)$ with TRG which is defined on the space of 
the $H$-bundle over $M$.
\end{rema} 

\subsection{Foliations with the zero structural Lie algebra}
The following proposition is proved using Theorem~\ref{t6}
by an analogy with \cite[Proposition 7]{ZhR}.

\begin{prop}\label{p2} Let $(M, F)$ be a foliation with TRG admitting 
an Ehresmann connection. Suppose that $\frak g_0(M,F)=0.$ Let
$\pi_b\colon{\cal R}\to W$ be the basic fibration. Then:
\begin{itemize} \item[(i)] the map $\Phi^W\colon W\times H\to W\colon
(w,a)\mapsto\pi_{\rm b} (R_a(u))\;\; \forall(w,a)\in W\times
H,\;\;\\   \forall u\in\pi_{\rm b} ^{-1}(w)$ defines a smooth locally
free action of the Lie group $H$ on the basic manifold $W;$

\item[(ii)] there is a homeomorphism $s\colon M/F\to W/H$ between
the leaf space $M/F$ and the orbit space $W/H$ satisfying the
equality $q\circ \pi_b=s\circ r\circ\pi,$  where $q\colon W\to
W/H$ and $r\colon M\to M/F$ are the quotient maps; 
\item[(iii)] the equality
$\pi_b^*\omega^W=\omega$ defines an $\mathbb{R}^k$-valued
non-degenerate $1$-form $\omega^W$ on $W$ such that
$\omega^W(A^*_W)=A$, where $A^*_W$ is the fundamental vector
field on $W$ defined by an element $A\in\mathfrak h\subset\mathbb{R}^k.$
\end{itemize}
\end{prop}

\section{Rigid geometries on leaf manifolds }

We want to emphasize that Sections~\ref{ss4.1} and \ref{ss4.2} $(M, F)$ 
is any smooth foliation. We do not assume the existence of an Ehresmann connection
for $(M, F)$.

\subsection{Generalized manifolds of leaf spaces of foliations}\label{ss4.1}

Consider a smooth foliation $(M, F)$ of codimension $n$ on $(n+m)$-dimensional manifold $M$.
Denote by $r:M\to M/F$ the quotient map onto the leaf space.
Refferring to \cite{Los1} let us consider the category $\mathbb R_n$
with objects open submanifolds of $\mathbb R^n$, where morphisms $Hom (U,V)$ are
diffeomorphisms $f:U\to V$ onto $f(U)\subset V$.

At every point $x\in M$ there exists a chart $(V,\varphi)$ of $M$ adapted to $(M,F)$.
This means that $\varphi(V) = W\times U\subset\mathbb R^{m+n}\cong\mathbb R^{m}\times\mathbb R^{n}$,
where $W$ and $U$ are open subsets in $\mathbb R^{m}$ and $\mathbb R^{n}$ respectively
and $\varphi(x) = (y,z)\in W\times U$. Denote by $pr:\mathbb R^{m}\times\mathbb R^{n}\to\mathbb R^{n}$ 
the canonical projection. The fibers of the submersion $pr\circ\varphi: V\to U$ belong to leaves 
of the foliation $(M, F)$. 

Let $r: M\to M/F$ be the projection onto the leaf space $M/F$ of $(M,F)$. 
For a fixed $y\in W$ denote by $j: U\to W\times U$ the embedding such that 
$j(u)= (y,u)\in W\times U$ $\,\, \forall u\in U.$ The pair $(U,k)$, where 
$k:= r|_U\circ\varphi^{-1}\circ j: U\to M/F$, is called a 
$\mathbb R_n$-chart (or chart) on $M/F$.

\begin{defi}\label{d10} Two charts $(U', k')$ and $(U'', k'')$ on $M/F$ for which 
$k'(U')\cap k''(U'')\neq\emptyset$, are called compatible if for each point $z\in k'(U')\cap k''(U'')$
there exists a chart $(U,k)$, $z\in k(U)$, with two morphisms $h':U\to U'$
and $h'':U\to U''$ in the category $\mathbb R_n$ satisfying the following conditions
$k'\circ h' = k$ and $k''\circ h'' = k$.
\end{defi}

\begin{defi}\label{d11} A smooth atlas on $M/F$ is a family of charts 
${\mathcal A} =\{(U_i,k_i)\,|\,i\in J\}$ satisfying the following two conditions:

1) the set $\{k_i(U_i)\,|\, i\in J\}$ is a covering of $M/F$, i.e. $\cup_{i\in J}k_i(U_i) = M/F$;

2) every two charts from $\mathcal A$ are compatible. 

A smooth atlas ${\mathcal A}$ is maximal if it is maximal relatively inclusion.
\end{defi}

\begin{defi}\label{d12} A pair $(M/F, {\mathcal A})$, where $\mathcal A$ is a maximal
atlas on $M/F$ is referred to as a leaf manifold. This leaf manifold is called 
induced by the foliation $(M, F)$ and it is denoted by $\mathcal N$.
The number $n$ is called the dimension of $\mathcal N$.
\end{defi}

Any atlas $\mathcal A$ defines the maximal atlas $\widehat{\mathcal A}$ as the set of charts which
are compatible with all charts from $\mathcal A$, hence $\mathcal A$ defines a smooth structure
on $M/F$.

\begin{defi}\label{d13} Let $({\mathcal N_1},{\mathcal A_1})$ and $({\mathcal N_2},{\mathcal A_2})$ 
be two $n$-dimensional leaf manifolds. A morphism of $\mathcal N_1$ to $\mathcal N_2$ is a map
$h:{\mathcal N}_1\to {\mathcal N}_2$ such that for each chart 
$(U,k)\in {\mathcal A}_1$ the pair $(U,h\circ k)$ is a chart from ${\mathcal A}_2.$ 
\end{defi}

Denote by $\mathfrak F$ the category of leaf manifolds of a fixed dimension $n$. 

\begin{defi}\label{d14} A foliation $(M', F')$ is referred to as associated with a leaf manifold
$\mathcal N$ if $(M', F')$ induces a leaf manifold $\mathcal N'$ on the leaf space $M'/F'$,
and $\mathcal N'$ is isomorphic to $\mathcal N$ in the category $\mathfrak F$.
\end{defi}
\begin{rema}
We emphasize that Definition~\ref{d12} is equivalent to the definition a smooth structure
on a leaf space in the sense of \cite{Los}.
It is known \cite[Theo\-rem~2]{Los} that for every foliation there exists a smooth atlas on the leaf space.
\end{rema}

\begin{rema}
On a leaf manifold $\mathcal N$ the tangent bundle $T\mathcal N$, differentials and codifferentials
of smooth maps, vector-valued forms are defined \cite{Los1}, \cite{Los2}. 
\end{rema}

\subsection{The pseudogroup approach to leaf manifolds}\label{ss4.2}
Let $(\mathcal N, \mathcal A)$ be a leaf manifold and $(M, F)$ be an associated foliation. 
Then the topological space of $\mathcal N$ is the leaf space $M/F$. 
Denote by ${\mathcal H} = {\mathcal H}(M,F)$ the holonomy pseudogroup of $(M, F)$. 
If the foliation $(M, F)$ is given by an $T$-cocycle $\{U_i, f_i, \{\gamma_{ij}\}\}_{i,j\in J}$, 
then ${\mathcal H}$ is generated by local diffeomorphisms $\gamma_{ij}$, and 
$T = \cup_{i\in J} f_i(U_{i})$. Let $T/\mathcal H$ be the quotient space. with the 
quotient map $q: T\to T/\mathcal H$. It is easy to show that $q: T\to T/\mathcal H$ is
continuous and open map. There exists a homeomorphism $\theta: M/F\to T/\mathcal H$
defined by the equality $\theta([L]):= [{\mathcal H}.v]$, where $[L]$ is a leaf $L = L(x)$
considered as a point of $M/F$ and $[{\mathcal H}.v]\in T/\mathcal H$ is the orbit 
${\mathcal H}.v$ of a point $v = f_i(x)$ for some submersion $f_i$ from the $T$-cocycle. 
Let us identify through $\theta$ the topological spaces $M/F$ and $T/\mathcal H.$

Consider a chart $(V,\psi)$ of the manifold $T$. Then $U=\psi(V)$ is an open subset in $\mathbb{R}^n.$
It is easy to see that the set $\mathcal B$ formed by the charts $(U,k)$ where $U=\psi(V)$ and
$k = q\circ\psi^{-1}: U\to M/F$ is an atlas on $M/F$, and $\mathcal B$ is compatible with
the atlas $\mathcal A$ of $\mathcal N$ induced by the foliation $(M, F)$.

Now let us show that the atlas $\mathcal B = \{(U_i,k_i)\,|\, i\in J\}$ defines 
a pseudogroup $\widetilde{\mathcal H}$ of local diffeomorpphisms of the $n$-dimensional manifold 
$\widetilde{T} = \coprod_{i\in J}U_i$. Let $(U_i, k_i)$ and $(U_j, k_j)$ be two charts from 
$\mathcal A$ and $z\in k_i(U_i)\cap k_j(U_j)$. According to the compatibility of these charts 
there exist a chart $(U,k)$ such that $z\in k(U)$ and two morphisms 
$h_i: U\to U_i$, $h_j: U\to U_j$ in the category $\mathbb R_n$ satisfying 
the equalities $k_i\circ h_i = k$ and $k_j\circ h_j = k$. Note that $h_i(U)$ 
and $h_j(U)$ are open subsets of $\widetilde{T}$. Therefore the map
$\widetilde{\gamma}_{ij}:= h_i\circ h_j^{-1}|_{h_j(U)}: h_j(U)\to h_i(U)$ is a local 
diffeomorphism of $\widetilde{T}$. The set $\{\widetilde{\gamma}_{ij}\,|\, i,j\in J\}$ 
generates a pseudogroup $\widetilde{\mathcal H}$ of local diffeomorphisms of $\widetilde{T}$. 
This pseudogroup is equivalent to the holonomy pseudogroup of any foliation
associated with the leaf manifold $\mathcal N$ in the terms of the work \cite{Sal}. 
Since the holonomy pseudogroup of a foliation is defined up to an equivalence, 
we have the following statement.
 
\begin{prop}\label{p3} Any leaf manifold $\mathcal N$ with an associated foliation $(M, F)$ 
may be defined by the indicated above atlas $\mathcal B$ defined by the
holonomy pseudogroup ${\mathcal H}(M,F)$ of $(M, F)$.

Every two foliations $(M, F)$ and $(M', F')$ associated with  $\mathcal N$  
have the same holonomy pseudogroup.
\end{prop}

Thus without loss of generality we consider the pseudogroup on $\widetilde{T}$ defined by the atlas
$\mathcal B$ as the holonomy pseudogroup of an associated foliation $(M, F)$ and use notations 
$\widetilde{T} = T$, $\widetilde{\mathcal H} = {\mathcal H} = 
{\mathcal H}({\mathcal N}) = {\mathcal H}(M, F).$

In the case when the smoothness is $C^\infty$, we may define the algebra of $C^\infty$-smooth 
functions $C^\infty(\mathcal N)$ on $\mathcal N$ as the algebra of $\mathcal H$-invariant 
$C^\infty$-smooth functions on $T.$ The Lie algebra of $C^\infty$-smooth
vector fields $\mathfrak{X}(\mathcal N)$ on $\mathcal N$ is defined as the Lie algebra of all 
derivations of the algebra of functions $C^\infty(\mathcal N)$. 

\subsection{Rigid geometries on leaf manifolds and their structural Lie algebras}\label{ss4.3}

\paragraph{The proof of Theorem~\ref{t1}}
Consider any $n$-dimensional leaf manifold $\mathcal N$. Let $(M, F)$ be 
an associated foliation, then $M/F$ is the topological space of $\mathcal N$. Assume that 
$(M, F)$ admits a transverse rigid geometry $\xi = (P(T,H),\beta)$. This is equivalent 
to the existence of a lifted foliation $({\mathcal R},{\mathcal F})$, where 
${\mathcal R}(M,H)$ is a principal $H$-bundle with the projection 
$\pi:{\mathcal R}\to M$, and the $\mathbb R^k$-valued $1$-form $\omega$ on 
$\mathcal R$ satisfying the conditions $(i)-(iii)$ from Section~\ref{ss3.3}. Denote by 
${\mathcal R}_{\mathcal F}$ the leaf manifold induced by $({\mathcal R},{\mathcal F})$.

A leaf $L$ considered as a point of the leaf space is denoted by $[L].$ Let $r:M\to\mathcal N$ 
and $r_{\mathcal F}: {\mathcal R}\to\mathcal R_{\mathcal F}$ be the projections onto the leaf manifolds 
$\mathcal N$ and $\mathcal R_{\mathcal F}$ respectively. 
Since $\pi: {\mathcal R}\to M$ is a morphism in the foliation category, 
the following map $\pi_{\mathcal F}: \mathcal R_{\mathcal F}\to{\mathcal N}: [\mathcal L]\mapsto
[\pi(\mathcal L)], \mathcal L\in\mathcal F$, is defined and satisfies the equality
$\pi_{\mathcal F}\circ r_{\mathcal F} = r\circ\pi.$ 

Due to $H$-invariance of $({\mathcal R},{\mathcal F})$, the map
$$\Phi:  {\mathcal R}_{\mathcal F}\times H\to{\mathcal R}_{\mathcal F}: [\mathcal L]
\mapsto[R_a(\mathcal L)],\, {{\mathcal L}}\in{\mathcal F},\,\,\, a\in H,$$
defines a right action of the Lie group $H$ on $\mathcal R_{\mathcal F}$. The isotropy subgroup
$H_{[\mathcal L]} = \{a\in H\,|\, R_a({\mathcal L}) = \mathcal L\}$ is a discrete subgroup of $H$
as the deck transformation group of the regular covering $\pi_{\mathcal L}: {\mathcal L}\to L.$
Therefore the action $\Phi$ of $H$ on $\mathcal  R_{\mathcal F}$ is locally free,
i.e. all isotropy groups are discrete subgroups of the Lie group $H$.
The orbit space ${\mathcal R}_{\mathcal F}/H$ is homeomorphic to $\mathcal N$. Let us identify 
$\mathcal N$ with ${\mathcal R}_{\mathcal F}/H$. We use the notation 
${\mathcal R}_{\mathcal F}({\mathcal N},H)$ for the quotient map 
$\pi_{\mathcal F}:{\mathcal R}_{\mathcal F}\to{\mathcal N}\cong{\mathcal R}_{\mathcal F}/H.$

Let $(M, F)$ and $(M', F')$ be two associated foliations with the same
leaf manifold $\mathcal N$. According to Proposition~\ref{p3} they have the common holonomy
pseudogroup ${\mathcal H} = \{\gamma_{ij}\,|\,i,j\in J\}$ of local automorphisms of the
rigid geometry $\xi = (P(T,H),\beta)$. Due to the efficiency of the
rigid geometry $\xi$ for each $\gamma_{ij}$ there exists a unique local automorphism $\Gamma_{ij}$ 
of $\xi$ lying over $\gamma_{ij}$ relative to the projection $p: P\to T.$ Emphasize that 
${\mathcal S} = \{\Gamma_{ij}\,|\, i, j\in J\}$ is the pseudogroup of the both lifted foliations
$(\mathcal R,\mathcal F)$ for $(M, F)$ and $(\mathcal R',\mathcal F')$ for $(M', F')$. 
So we will use the following notation ${\mathcal S} = \mathcal S({\mathcal N},\xi).$ Since 
$\mathcal R_{\mathcal F} = P/\mathcal S$ is defined by $\mathcal S$, the leaf manifold 
$\mathcal R_{\mathcal F}$ is not depend on the choice of an associated foliation. 

Every $\Gamma_{ij}$ acts freely on $P$ as a local automorphism of the parallelizable
manifold $(P,\beta)$. The free action of the pseudogroup 
${\mathcal S} = \{\Gamma_{ij}\,|\, i, j\in J\}$ on $P$ implies the free action on $TP$ of 
the pseudogroup ${\mathcal S}_*: = \{\Gamma_{ij*}\,|\, i, j\in J\}$
formed by differentials of local transformations belonging to $\mathcal S,$ and 
$T{\mathcal R}_{\mathcal F} = TP/{\mathcal S_*}.$ Therefore 
the tangent space $T_z\mathcal R_{\mathcal F}$ is a vector space, and 
the dimension of $T_z\mathcal R_{\mathcal F}$ is equal to $k = \dim(P)$ at any point 
$z\in\mathcal R_{\mathcal F}.$ As $\Gamma_{ij}^*\beta = \beta,\,\, i, j\in J,$ a 
non-degenerate $\mathbb R^k$-valued $1$-form $\alpha$ is defined on 
$\mathcal R_{\mathcal F}$ and satisfies the equality $\mu^*\alpha =\beta$, where 
$\mu: P\to P/{\mathcal S} = \mathcal R_{\mathcal F}$ is the quotient map.
We emphasize that $\alpha$ coincides with the $1$-form defined by the following
equality $\pi_{\mathcal F}^{*} \alpha = \omega$, where $\omega$ is the basic 
$\mathbb R^k$-valuated $1$-form on $\cal R $ satisfying the conditions 
$(i) - (iii)$ in Section~\ref{ss3.3}. 

Since ${(\mathcal R,\mathcal F)}$ is a Riemannian foliation with the Ehresmann connection 
$\widetilde{\mathfrak M} = \pi^{*}\mathfrak M$ where $\mathfrak M$ is an Ehresmann connection 
for $(M, F)$, according to \cite[Proposition 2]{ZhG} the holonomy pseudogroup $\mathcal S$
of ${(\mathcal R,\mathcal F)}$ is complete. By \cite[Theorem 3.1]{Sal} the structural Lie algebra
$\mathfrak g_0(\mathcal S)$ is defined, and this Lie algebra is isomorphic to the structural Lie algebra 
$\mathfrak g_0(\mathcal R, \mathcal F)$ of $(\mathcal R, \mathcal F)$. According to Definition \ref{d8}
$\mathfrak g_0(M, F): = \mathfrak g_0(\mathcal R, \mathcal F)$. As ${\mathcal S} = \mathcal S({\mathcal N},\xi)$, 
then $\mathfrak g_0=\mathfrak g_0(\mathcal S)$ is not depend on the choice of an associated foliation $(M, F)$. 
$\square$

\begin{defi}\label{d:15} The pair $\zeta = ({\mathcal R}_{\mathcal F}({\mathcal N},H),\alpha)$
defined in the proof of Theorem \ref{t1} is called the rigid geometry on the leaf manifold 
$\mathcal N$ modelled on the transverse rigid geometry $\xi = (P(N,H),\beta)$ of an associated 
foliation $(M, F)$. The structural Lie algebra $\mathfrak g_0$ of $(M, F)$
is called by the structural Lie algebra of the rigid geometry $\zeta$ on $\mathcal N$ and is denoted by
$\mathfrak g_0 =\mathfrak g_0(\zeta).$
\end{defi}
According to the proof of Theorem~\ref{t1} Definition~\ref{d:15} is correct, i.e. $\mathfrak g_0$
does not depend on the choice of the associated foliation $(M, F)$.

We want to emphasize that the rigid geometry 
$\zeta = ({\mathcal R}_{\mathcal F}({\mathcal N},H),\alpha)$ and the structural Lie algebra 
$\mathfrak g_0$ depend only on the pseudogroup $\mathcal H$ of local automorphisms of the 
transverse geometry $(T,\xi)$ on which $(M, F)$ is modelled.

\section{Automorphisms of rigid geometries on leaf manifolds}
\subsection{The category of rigid geometries on leaf manifolds}

\begin{defi}\label{d:16} Let $\zeta = ({\mathcal R}_{\mathcal F}({\mathcal N},H),\alpha)$
and $\zeta' = ({\mathcal R}'_{\mathcal F'}({\mathcal N}',H'),\alpha')$
be two rigid geometries on $n$-dimensional leaf manifolds $\mathcal N$ and $\mathcal N'$ respectively.
A map $\Gamma: {\mathcal R}_{\mathcal F}\to{\mathcal R}'_{\mathcal F'}$ is called a morphism 
$\zeta\to\zeta'$ if $H = H'$ and $\Gamma$ satisfies the following two conditions:
$$1)\, \Gamma^*\alpha' =\alpha,\,\,\,\, 2) \,\Gamma\circ R'_a = R_a\circ\Gamma\,\,\, \forall a\in H. $$
\end{defi}
 Thus rigid geometries on leaf manifolds form the category $\mathfrak R\mathfrak F$.
Let $\mathfrak R\mathfrak F_0$ be the full subcategory of $\mathfrak R\mathfrak F$ objects of which 
have zero structural Lie algebra.

Every $\Gamma\in Mor(\zeta,\zeta')$ defines the projection $\gamma:{\mathcal N}\to {\mathcal N}'$ 
such that $\pi'_{\mathcal F'}\circ\Gamma =\gamma\circ\pi_{\mathcal F}$, where 
$\pi_{\mathcal F}: {\mathcal R}_{\mathcal F}\to {\mathcal N}\cong{\mathcal R}_{\mathcal F}/H$ and 
$\pi'_{\mathcal F'}: {\mathcal R'}_{\mathcal F'}\to {\mathcal N'}\cong{\mathcal R'}_{\mathcal F'}/H$
are the quotient maps. There is a covariant functor
${\mathcal K}:\mathfrak R\mathfrak F\to\mathfrak F$ forgetting rigid geometries. Hence
${\mathcal K}(\zeta) = {\mathcal N}$ for each $\zeta\in Ob(\mathfrak R\mathfrak F)$ and
${\mathcal K}(\Gamma) = \gamma$ for every $\Gamma\in Mor(\zeta,\zeta')$.
Let $\mathfrak F_0: = {\mathcal K}(\mathfrak R\mathfrak F_0).$

\subsection{Proof of Theorem~\ref{t2}.}
Consider a rigid geometry 
$\zeta = ({\mathcal R}_{\mathcal F}({\mathcal N},H),\alpha)\in\mathfrak R\mathfrak F_0$ 
on $n$-dimensional leaf manifold $\mathcal N$, then $\mathfrak g_0(\zeta) = 0.$ 
In accordance with Proposition~\ref{p2} in this case the leaves of the lifted foliation
$({\mathcal R},{\mathcal F})$ are fibers of the locally trivial basic fibration 
$\pi_b: {\mathcal R}\to W$, and the leaf space ${\mathcal R}/{\mathcal F}$ is $W$. Therefore 
there exists an isomorphism $f: {\mathcal R}_{\mathcal F}\to W$ in the category $\mathfrak F_0$ of
the leaf manifold ${\mathcal R}_{\mathcal F}$ onto the manifold $W$ of dimension $\dim(W) = n+s$, 
where $n = \dim(\mathcal N)$ and $s = \dim(H)$. Moreover, $f^*\omega^W = \alpha$, hence $f$ is an isomorphism
of parallelelizable manifolds $({\mathcal R}_{\mathcal F},\alpha)$ and $(W,\omega^W)$.
Let us identify $({\mathcal R}_{\mathcal F},\alpha)$ with $(W,\omega^W)$ through $f$,
then ${\mathcal R}_{\mathcal F} = W$ and $\alpha = \omega^W$. Therefore the statement~$(i)$ 
of Theorem~\ref{t3} follows from Propositions~\ref{p2}.

As it is well known, the automorphism group 
$A(W,\alpha) =\{h\in\text{\rm Diff}(W)\,|\, f^*\alpha =\alpha\}$
of the parallelisable manifold $(W,\alpha)$ admits a Lie group structure, 
and its dimension is not grater than $\dim(W).$
According to Definition~\ref{d:16}, the group $Aut(\zeta)$ 
of all automorhisms of $\zeta$ in the category $\mathfrak F$ is equal to 
$$Aut(\zeta) = \{h\in Diff(W)\,|\, h^*\alpha=\alpha, \,\,
R^W_a\circ h = h\circ R^W_a \,\, \forall a\in H\}.$$ 
Hence, $Aut(\zeta) = \{h\in A(W,\alpha)\,|\,R^W_a\circ h = h\circ R^W_a, a\in H\}.$
This implies that $Aut(\zeta)$ is a closed subgroup of the Lie group
$A(W,\alpha)$. It means that $Aut(\xi)$ admits a Lie group structure,
and its dimension is not grater than $\dim(W).$ Since $Aut(\zeta)$ is
a transformation group, $Aut(\zeta)$ is equipped with the compact-open topo\-logy 
and a Lie group structure on $Aut(\zeta)$ is unique
\cite[Theorem~VI]{Pal}. Thus statements~$(ii)$ and $(iii)$ 
of Theorem~\ref{t2} are proved. \,\,\,
$\square$

\subsection{Proof of Theorem~\ref{t3}.}

A smooth foliation $(M, F)$ is called proper if
each its leaf is an embedded submanifold in $M$. A subset $M_0$ of $M$ is
said to be saturated if it is a union of leaves of $(M, F)$.

Recall that a topological space satisfies the separation axiom $T_0$ if for any
different points $a$ and $b$ there exists a neighborhood at least one of them contains no 
other point. As it is known \cite[Lemma 4.1]{ZhT}, a foliation is proper if and only if its leaf 
space satisfies the separation axiom $T_0$.
Assume that a leaf manifold $\mathcal N$ satisfies the separation axiom $T_0$.
Hence according to the mentioned above lemma every associated foliation $(M, F)$ is proper. 
Assume that $(M, F)$ is a foliation with transversal rigid geometry $\xi$. According 
to Theorem~\ref{t1} the induced rigid geometry $\zeta$ on $\mathcal N$ is defined. 
Consider the lifted foliation $(\mathcal R,\mathcal F)$ for $(M, F)$. 
By the assumption $(M, F)$ admits an Ehresmann connection $\mathfrak M$. 
Preserving the notations used above we denote by $\pi:{\mathcal R}\to M$ 
the projection of the $H$-bundle ${\mathcal R}(M,H)$. Then 
$\widetilde{\mathfrak M} = \pi^*\mathfrak M$ is an Ehresmann connection for
$(\mathcal R,\mathcal F)$. 

It is known that any foliation has a leaf with the trivial germ holonomy group.
Thus there is a proper leaf with the trivial germ holonomy group of the foliation 
$(M, F)$. This implies the existence of a proper leaf of the $e$-foliation  
$(\mathcal R,\mathcal F)$ admitting an Ehresmann connection. 
Therefore $(\mathcal R,\mathcal F)$ is also a proper foliation. Since the closure 
$\overline{\mathcal L}$ of a leaf ${\mathcal L}\in{\mathcal F}$ is a minimal set,
in the case when $\overline{\mathcal L}\neq\mathcal L$, the closure $\overline{\mathcal L}$ 
contains only non-proper leaves. Hence it is necessary $\overline{\mathcal L} = \mathcal L.$ 
This implies that the structural Lie algebra $\mathfrak g_0 = \mathfrak g_0(\zeta)$ is zero. 
Therefore the rigid geometry $(\mathcal N,\zeta)$ satisfies Theorem~\ref{t2} and 
the statement $1)$ of Theorem~\ref{t3} is proved.

Similarly to \cite[Theorem 1.1]{ZhT} for proper Cartan foliations admitting Ehresmann conections,
we prove the following statement. 
\begin{theo}\label{t7} Let $(M, F)$ be an arbitrary proper foliation of codimension $n$ 
with transverse rigid geometry $\xi = (P(T,H),\beta)$ admitting an Ehresmann connection. 

Then there exists a not necessarily connected, saturated, dense open subset $M_0$ of $M$
such that the induced foliation $(M_0, F|_{M_0})$ is formed
by the fibers of a locally trivial fibration $p : M_0\to B$ with the standard fiber $L_0$ 
over a smooth $n$-dimensional (not necessarily Hausdorff) manifold $B.$
\end{theo}

Now the statement $2)$ of Theorem~\ref{t3} follows from Theorem~\ref{t7}.
\,\,$\square$

\end{document}